\documentclass[twoside]{article}
\usepackage[dvips]{graphicx}
\usepackage{color}
\usepackage{enumerate}% permet de choisir la forme des numeros ds \enumerate
\usepackage{amssymb}
\usepackage{multirow}
\usepackage[colorlinks=false,urlcolor=blue,citecolor=blue,linkcolor=blue,bookmarks=true, pdfstartview=FitH, bookmarksopen=true]{hyperref}

\usepackage{a4}

\usepackage{theorem,amssymb,graphicx}

\newtheorem{theorem}{ Theorem}[section]

\newtheorem{definition}{Definition}[section]
\newtheorem{example}{Example}[section]
\newtheorem{lemma}{Lemma}[section]
\newtheorem{proposition}{Proposition}[section]
\newtheorem{remark}{Remark}[section]

\newcommand{\qed}{\hskip 10pt $\Box$}
\newenvironment{proof}{\par {\it Proof.\hskip 5pt}}{\hfill \qed \par}
\newenvironment{proofof}[1]{\par {\it Proof of #1.} \hskip 5pt}{\hfill\qed\par}

\def\NN{\mathbb N}
\def\R{\mathbb R}
\def\1{\mathrm{1\hspace{-2pt}I}}
\def\rma{{\mathbb R}_{\max}}
\def\rhom{{\rho}_{\max}}
\def\mp{max-plus }
\def\mE{\mathcal E}
\def\M{\mathcal M}
\def\mP{\mathcal{P}}
\newcommand{\junk}[1]{}%
\def\G{\mathcal{G}}
\def\Gc{\mathcal{G}^c}
\def\P{\mathbb{P}}

\title{Memory loss property for products of random matrices in the max-plus  algebra} % insert title - use \\ if it requires 
\author{G. MERLET}

\begin{document}
\maketitle
\begin{abstract}
Products of random matrices in the max-plus  algebra are used as a model for a class of discrete event dynamical systems.  This can model a wide range of systems including train or queuing networks, job-shop, timed digital circuits or parallel processing systems.

Some stability results have been proved under the so-called memory loss property.

When the random matrices are i.i.d, we prove that the memory loss property is generic in the following sense: if it is not fulfilled, the support of the common law of the random matrices is included in a union of finitely many affine hyperplanes and in the discrete case the atoms of the measure are linearly related.
\end{abstract}

\normalsize

\section{Introduction}
This paper will focus on sequences of $\R^k$ valued random variables $\left(X_n\right)_{n\in\NN}$ defined by a stationary sequence  of random matrices $\left(A(n)\right)_{n\in_\NN}$, an initial condition $X_0$ and the recurrence relation
\begin{equation}\label{EqRecCoeffs}
(X_{n+1})_i= \max_j\left( A_{ij}(n)+(X_n)_j\right).
\end{equation}

Some stability results (J.~Mairesse~\cite{Mairesse}, S.~Gaubert and D.~Hong~\cite{GaubertHong}, G.~M.~\cite{limtop, TclGM}) for such sequences have been proved under the so-called memory loss property (MLP).

When the matrices are i.i.d., we prove that the property is generic in the following sense: if the sequence has not the MLP, then the support of the common law of the random matrices is included in a union of finitely many affine hyperplanes and in the discrete case the atoms of the measure are linearly related.

This article is divided into four sections. In the first one, we introduce the \mp algebra and illustrate its modeling power with one simple example and many references.
In the second section, we define the MLP and explain its interest, then we state and comment our main results. In the third section, we remind the asymptotic theory of matrices in the \mp algebra and we prove a needed extension of those results. In the last section, we prove the main results.

\section{Modeling}
\subsection{Practical situations}
This class of system can model a wide range of situations.  Among other examples it has been applied to queuing networks (J. Mairesse~\cite{Mairesse}, B.~Heidergott~\cite{CaractMpQueuNet}), train networks (B.~Heidergott and R.~De~Vries~\cite{HeidergottDeVriesPubTransNet}, H.~Braker~\cite{braker}) or job-shop (G. Cohen et al.~\cite{cohen85a}). It also computes the daters of some task resources models (S.~Gaubert and J.~Mairesse~\cite{gaumair95}) and  timed Petri nets including events graphs (F.~Baccelli~\cite{Baccelli}) and 1-bounded Petri nets (S.~Gaubert and J.~Mairesse~\cite{GaubertMairesseIEEE}). The role of the max operation is synchronizing different events. For developments on the max-plus modeling power, see the books by F.~Baccelli et al.~\cite{BCOQ} or B.~Heidergott et al.~\cite{MpAtWork}.\\

To illustrate our results, below is a very simple system ruled by \mp~equations. 

\begin{example}\label{exprod}
The process consists in two tasks, performed on the same piece in two different places. Before the first task, the piece is put on a kart. It is removed of it after completion of the second task. There are two karts, used one after the other.

We denote by $(X_n)_1$ the date of completion of the first task for the $n^\mathrm{th}$ time and by $(X_n)_2$ the date of completion of the second task for the $(n-1)^\mathrm{th}$ time. For any~$i$, $\gamma_i(n)$ is the duration of the $i^\mathrm{th}$ task for the $n^\mathrm{th}$ time.  The kart takes time~$t_1(n)$ to go from first to second place for the $n^\mathrm{th}$ time and time~$t_2(n)$ to go back to the first place. This is sum up in figure~\ref{BuildSys}.
\begin{figure}[htbp]
\begin{center}
\caption{A simple production system}\label{BuildSys}\ \\
\begin{picture}(0,0)%
\includegraphics{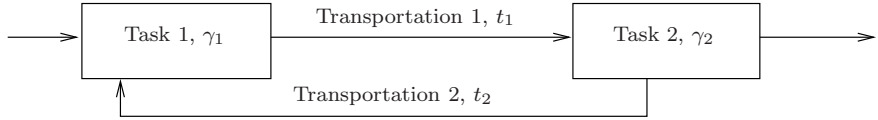}%
\end{picture}%
\setlength{\unitlength}{3947sp}%
\begingroup\makeatletter\ifx\SetFigFont\undefined%
\gdef\SetFigFont#1#2#3#4#5{%
  \reset@font\fontsize{#1}{#2pt}%
  \fontfamily{#3}\fontseries{#4}\fontshape{#5}%
  \selectfont}%
\fi\endgroup%
\begin{picture}(5480,733)(698,-827)
\put(1442,-342){\makebox(0,0)[lb]{\smash{{\SetFigFont{8}{9.6}{\rmdefault}{\mddefault}{\updefault}Task 1, $\gamma_1$}}}}
\put(4513,-342){\makebox(0,0)[lb]{\smash{{\SetFigFont{8}{9.6}{\rmdefault}{\mddefault}{\updefault}Task 2, $\gamma_2$}}}}
\put(2647,-224){\makebox(0,0)[lb]{\smash{{\SetFigFont{8}{9.6}{\rmdefault}{\mddefault}{\updefault}Transportation 1, $t_1$}}}}
\put(2505,-720){\makebox(0,0)[lb]{\smash{{\SetFigFont{8}{9.6}{\rmdefault}{\mddefault}{\updefault}Transportation 2, $t_2$}}}}
\end{picture}%
\end{center}
\end{figure}

A task is performed as soon as the place is free and there is an available kart. This means that the sequence $\left(X_n\right)_{n\in\NN}$ is ruled by equations
$$(X_{n+1})_1-\gamma_1(n+1)=\max\left((X_n)_1,(X_n)_2+t_2(n-1)\right)\mathrm{ and }(X_{n+1})_2-\gamma_2(n)=\max\left((X_n)_1+t_1(n),(X_n)_2\right),$$
in which we recognize equation~(\ref{EqRecCoeffs}) with matrix 
$A(n)=\left(\begin{array}{cc}
\gamma_1(n+1)&\gamma_1(n+1)+t_2(n-1) \\
\gamma_2(n)+t_1(n)&\gamma_2(n) \end{array}\right).$

It is natural to assume that the sequence $\left(t_1(n),t_2(n),\gamma_1(n),\gamma_2(n)\right)_{n\in\NN}$ is stationary. Therefore, so is $\left(A(n)\right)_{n\in\NN}$ . From now on, we also assume that the '$A(n)$'s are independent. It is the case if the variations in the durations are caused by phenomena completely external to the process.
\end{example}

\subsection{Max-plus algebra}
Sequences satisfying equation~(\ref{EqRecCoeffs}) are best understood by introducing the so-called \mp algebra, which is actually a semiring.
\begin{definition}
The \mp semiring $\rma$ is the set $\R\cup\{-\infty\}$,  with  the max as a sum (i.e. 
$a \oplus b = \max(a, b)$) and the usual sum as a product (i.e. $a\otimes b = a + b$). In this semiring, the neutral elements are $-\infty$ and $0$.
\end{definition}

We also use the matrix and vector operations induced by the semiring structure.
For matrices $A, B$ with appropriate sizes, $(A\oplus B)_{ij} = A_{ij} \oplus B_{ij} = \max(A_{ij} , B_{ij})$, $(A\otimes B)_{ij} = 
\bigotimes_k A_{ik}\otimes B_{kj} = \max_k (A_{ik} + B_{kj} )$, and for a scalar $a$, 
$(a\otimes  A)_{ij} = a\otimes  A_{ij} = a + A_{ij}$. Given an integer $n$, we denote by $[1,n]$ the set $\{1,\cdots,n\}$.

\begin{definition}
Let $\mathcal{A}=\left(A(n)\right)_{n\in_\NN}$ be a stationary sequence of random ${k\times k}$-matrices with entries in $\rma$ and no row of $-\infty$.
 
We investigate the behavior of the sequence $\left(x^\mathcal{A}(n,x_0)\right)_{n\in_\NN}$ defined by $x^\mathcal{A}(0,x_0)=x_0$ and $$x^\mathcal{A}(n+1,x_0)=A(n)\otimes  x^\mathcal{A}(n,x_0).$$
\end{definition}

With those definitions the random variable $X_n$ of  equation~(\ref{EqRecCoeffs}) is exactly $x^\mathcal{A}(n,X_0)$.

\section{Memory loss property}
\subsection{Definition}
Let us remind a few basic facts about matrices with entries in $\rma$. All these facts can be checked by direct computations.
\begin{enumerate}[i)]
\item A matrix $A\in\rma^{k\times l}$ defines a max-plus linear map $\tilde{A}$ from $\rma^l$ to $\rma^k$ given by:
$$\forall x \in\rma^l, \tilde{A}(x):=A\otimes x,$$
and the product of the matrices corresponds to the composition of the maps.
If $A$ has no row of $-\infty$, it also defines a map from $\R^l$ to $\R^k$.
\item The image of $\tilde{A}$ is stable under the operations of $\rma^k$. It is  the max-plus moduloid spanned by the columns vectors of $A$. Like in usual linear algebra, we have $\forall x \in\rma^k,A\otimes x=\bigoplus_{j\in [1,k]}x_j\otimes A_{.j}\,. $
\end{enumerate}

\begin{definition}\label{defMLP}\ 
\begin{itemize}
\item A matrix $A$ has rank 1 if all its columns are proportional in the \mp sense and not all equal to~$(-\infty)^k$. It happens if and only if there are $a$ and $b$ in $\rma^k\backslash\{(-\infty)^k\}$ such that $\forall (i,j)\in [1,k]^2$, $A_{ij}=a_i+b_j$. We denote it by $rk(A)=1$.
\item A sequence  $\left(A(n)\right)_{n\in_\NN}$ is said to have the memory loss property (MLP) if there exist an $N\in\NN$ such that $\mathbb{P}\left( rk\left(A(N)\otimes  \cdots \otimes  A(1)\right)=1\right)>0$.
\end{itemize}
\end{definition}

The expression ``memory loss property'' was introduced by S.~Gaubert and D.~Hong~\cite{GaubertHong} for \mp automata. It expresses that $rk\left(A(n)\otimes  \cdots \otimes  A(1)\right)=1$ if and only if, for every pair $(i,j)$, the difference $x_i^\mathcal{A}(n,x_0)-x_j^\mathcal{A}(n,x_0)$ does not depend on the initial condition $x_0$. 

Let us review the three type of stability results that were proved under the MLP hypothesis.

J.~Mairesse introduced the MLP to show the result below:
\begin{theorem}\label{scd}(Mairesse~\cite[6.15]{Mairesse})
If $\mathcal{A}$ has the MLP, then for every $i,j\in [1,k]$ the sequences $\left(x_i^\mathcal{A}(n,x_0)-x_j^\mathcal{A}(n,x_0)\right)_{n\in\NN}$ and $\left(x_i^\mathcal{A}(n+1,x_0)-x_j^\mathcal{A}(n,x_0)\right)_{n\in\NN}$ converge in total variation, uniformly in~$x_0$.
\end{theorem}

It is known since J.E. Cohen~\cite{Cohen}, that $\left(x^\mathcal{A}_i(n,x_0)\right)_{n\in\NN}$ satisfies a law of large number, at least when the entries are finite and integrable. The limit of $\left(\frac{1}{n}\max_ix^\mathcal{A}_i(n,0)\right)_{n\in\NN}$ exists under the assumption that $\max_{i,j\in [1,k]}\left(A_{ij}(0)\right)^+\in\mathbb{L}^1$ (see~J.-M.~Vincent~\cite{vincent}) and  is called the Lyapunov exponent of $\left(A(n)\right)_{n\in\NN}$. The other two results show that the MLP ensure a kind of stability of the Lyapunov exponent.

If $A(1)$ takes only finitely many values, then the MLP obviously only depends on these values, and not on the law of $A(1)$. If the sequence defined by the uniform distribution on $(A_1,\cdots, A_t)$ has the MLP, then  every sequence $\left(A(n)\right)_{n\in\NN}$ of i.i.d random variables, such that \mbox{$\forall i\in [1,t], \mathbb{P}(A(n)=A_i)>0$} also has the MLP. That being the case, we have the following result on the Lyapunov exponent:
\begin{theorem}(Gaubert and Hong~\cite[4.1]{GaubertHong})
Let  $(A_1\cdots A_t)$ be a $t$-uple of matrices in $\rma^{k\times k}$ with no row of $-\infty$. We write $p>0$ for $\forall i\in [1,t], p_i>0$.
For every probability vector $p$, we call $L(p)$ the Lyapunov exponent of i.i.d random variables $A^p(n)$ such that 
$$\forall i\in [1,t], \mathbb{P}(A^p(1)=A_i)=p_i.$$
If there is one $p>0$ (or for every $p>0$), such that $\left(A^p(n)\right)_{n\in\NN}$  has the MLP, then $L$ is an analytical function of $p$. The analyticity domain only depends on the $t$-uple and  includes all $p>0$.
\end{theorem}

Eventually, I have proved~\cite{limtop, TclGM} that if $\mathcal{A}=\left(A(n)\right)_{n\in\NN}$ has the MLP, and the '$A_n$'s are sufficiently integrable and mixing, then $\left(x^\mathcal{A}(n,X_0)\right)_{n\in\NN}$ satisfies a central limit theorem. If the '$A_n$'s are independent, it also satisfies local limit theorem, renewal theorem and great deviation principle.

To conclude this section, we give an example of sequence without the MLP:
\begin{example}\label{ExNonMlp}
If there is a constant $c\in\R$, such that $\left(A(n)\right)_{n\in\NN}$ takes it values in the set 
$$A_c=\left\{\left.\left(\begin{array}{cc}\gamma_1&\gamma_1+t \\\gamma_2+t&\gamma_2 \end{array}\right)\right|\gamma_1\ge 0,\gamma_2\ge 0, \gamma_2 -\gamma_1=c,t\ge |c|\right\},$$ then the sequence $\left(A(n)\right)_{n\in\NN}$ has not the MLP.
\end{example}
\begin{proof}
Let us define vectors $u=(0,0)'$ and $v=(0,c)'$. Matrices in $A_c$  map $u$ on $(\gamma_1+t)\otimes v$ and $v$ on $(\gamma_2+t)\otimes u$.

Therefore, the image of $A(n)\otimes\cdots\otimes A(1)$ always contains $u$ and $v$, which are not proportional in the \mp~sense.
\end{proof}

\subsection{Statements of the results}
Since we want $x^\mathcal{A}$ to be finite, matrices will take their values in the set of $k\times k$ matrices with entries in $\rma$ and no row of $-\infty$, which we denote by $\mathcal{M}_k$. For convenience, we denote by $\mathcal{P}_k$ the set of the so-called primitive matrices, that is matrices $A\in\mathcal{M}_k$ such that there is an integer $n\in\NN$ such that $A^{\otimes n}$ has only finite entries.

We set a notion of linear forms and hyperplane on these space, which extend the usual notions on $\R^{k\times k}$.

\begin{definition}
For a finite set $Q$ and $\alpha\in\R^Q$, we define $f_\alpha:\rma^Q\rightarrow \rma$ by:
$$f_\alpha(V)=
\left\{\begin{array}{ll}
\sum_{\alpha_i\neq 0}\alpha_iV_i&\mathrm{ if }\forall i\in Q, \alpha_i\neq 0\Rightarrow V_i\in\R,\\
&\textrm{ where the operations are understood in the usual sense}\\
-\infty &\textrm{ otherwise.}
\end{array}\right. $$
Such a function is called a linear form on $\rma^Q$. A hyperplane (rep. affine hyperplane) of $\rma^Q$ is the set of zeros (resp.  the level line associated to a finite level) of a non-zero linear form.

A hyperplane of $\mathcal{M}_k$ (resp. $\mathcal{P}_k$) is the intersection of a hyperplane  of $\rma^{k\times k}$ with $\mathcal{M}_k$ (resp. $\mathcal{P}_k$).

A hyperplane of $\mathcal{M}_k\times \mathcal{P}_k $ is the intersection of a hyperplane  of $\left(\rma^{k\times k}\right)^2$ with $\mathcal{M}_k\times \mathcal{P}_k$.

When $Q=k\times k$, for every indices $i,j\in [1,k]$, we denote by $A^\circ_{ij}$ the function $A\mapsto A_{ij}$. 
That being the case, we write $f_\alpha=\sum \alpha_{ij}A^\circ_{ij}$ and this decomposition is unique.
\end{definition}

Our main results are the two following theorems:
\begin{theorem}\label{thgenedis}
For every $k\in \NN$, the set of pairs of matrices $(A,B)\in\mathcal{P}_k\times\mathcal{M}_k$ such that a sequence of i.i.d. random matrices that take values $A$ and $B$ with positive probability does not have the MLP is included in the union of finitely many hyperplanes of $\mathcal{P}_k\times\mathcal{M}_k$.

For such a pair $(A,B)$, a stationary sequence $\left(A(n)\right)_{n\in\NN}$ has the MLP, provided it satisfies the following relation:
$$\forall N\in\NN, \forall (A_i)_{i\in [1,N]}\in \{A,B\}^N, \P\left[\forall i\in  [1,N], A(i)=A_i\right]>0. $$
\end{theorem}

To have a similar result for continuous measures, we need to define the support of a measure. To this aim, we put the following distance on $\mathcal{M}_k$: $d(A,B)=\max\{\left.\left|\arctan(A_{ij})-\arctan(B_{ij})\right|\right|i,j\in [1,k] \}.$

\begin{theorem}\label{thgenesupp}
Let $\mu$ be a probability measure on $\mathcal{M}_k$ with support $S_\mu$. If $S_\mu\cap \mathcal{P}_k$ is not included in a union of finitely many affine hyperplanes of $\mathcal{P}_k$, then  a sequence of i.i.d. random matrices with law $\mu$ has the~MLP.

For such a probability  measure $\mu$, a stationary sequence $\left(A(n)\right)_{n\in\NN}$ has the MLP, provided  $A(1)$ has law $\mu$ and for every $N\in\NN$, the product set $(S_\mu\cap \mathcal{P}_k)^N$ is included in the support of $\left(A(n)\right)_{n\in [1,N]}$ 
\end{theorem}

\begin{example}[Application to example~\ref{exprod}]
Because of example~\ref{ExNonMlp}, it is obviously possible that the sequence ruling the system described in example~\ref{exprod} has not the MLP. It depends on the possible durations of the tasks and transportation times. However, Theorems~\ref{thgenesupp} and~\ref{thgenedis} imply that the sequence generically has the MLP: if it has not, the support of $\left(\gamma_1(2), \gamma_2(1), t_1(1), t_2(0)\right)$ is included in the union of finitely many affine hyperplanes and in the discrete case the atoms of the measure are linearly related.

Together with~\cite{limtop, TclGM}, it means that, if $\left(\gamma_1(2), \gamma_2(1), t_1(1), t_2(0)\right)$ is sufficiently integrable and not almost surely included in a union of finitely many affine hyperplanes, then $\left(X_n\right)_{n\in\NN}$ satisfy a central limit theorem, a local limit theorem, and a renewal theorem.
\end{example}

\subsection{Outlines of the proof}\ 

\noindent{\bf Notation}
We will denote  a $n$-uple of matrices by $(^i\hspace{-4pt}A)_{i\in [1,n]}$ instead of $(A_i)_{i\in [1,n]}$ to use indices for entries of matrices.

Theorem~\ref{thgenedis} (very partially) answers the following deterministic question: when does the semigroup generated by two matrices in $\rma^{k\times k}$ contains a matrix with rank~1? The asymptotic theory of matrices in $\rma^{k\times k}$ shows that some matrices, called \textbf{scs1-cyc1}, have powers of rank~1. Theorem~\ref{thgenedis} thus readily follows  from the next two lemmas which will be proved in sections~\ref{seclem} and~\ref{deter} respectively.
\begin{lemma}\label{graphegene}
For every pair $(A,B)\in \mathcal{P}_k\times\mathcal{M}_k$ outside a union of finitely many hyperplanes, there exist two integers $m$ and $n$ such that the matrix $A^{\otimes m}\otimes B\otimes A^{\otimes n}$ has only finite entries and its  critical graph has only one node. As a consequence, that matrix is \textbf{scs1-cyc1}.
\end{lemma}
\begin{lemma}\label{scs1-cyc1}
Every \textbf{scs1-cyc1} matrix with finite entries has a power with rank 1.
\end{lemma}

Theorem~\ref{thgenesupp} will be deduced from lemma~\ref{graphegene} and the following:
\begin{lemma}\label{puissrg1}
If $A$ is a matrix with finite entries whose critical graph has only one node, then there are a neighborhood $V$ of $A$ and an integer  $n$ such that: 
$$\forall (^i\hspace{-4pt}A)_{i\in [1,n]}\in V^n, rk(^1\hspace{-4pt}A\otimes\cdots\otimes ^n\hspace{-6pt}A)=1.$$
\end{lemma}

\begin{proofof}{theorem~\ref{thgenesupp}}\ \\
Every hyperplane $\alpha$ given by lemma~\ref{graphegene} is the kernel of a linear form $f^\alpha$ on $\mP_k\times\M_k$. This linear form can be written $f_1^\alpha+f_2^\alpha$ where $f^\alpha_1$ depends only on the first matrix and $f^\alpha_2$ on the second one.
Since $S_\mu\cap\mP_k$ is not included in $\cup_\alpha \ker f_1^\alpha$, there exists a matrix $A$ in $S_\mu\cap\mP_k$ such that $\forall \alpha, f_1^\alpha(A)\neq 0$. For every $\alpha$, the set $\{B\in\M_k|f_2^\alpha(B)=-f_1^\alpha(A)\}$ is an affine hyperplane  of $\M_k$ or the emptyset. Therefore, there exists a $B\in S_\mu$ such that $B\not\in\cup_\alpha\{B\in\R^{k\times k}|f_2^\alpha(B)=-f_1^\alpha(A)\}$. Eventually, $(A,B)\notin \bigcup_\alpha \ker f^\alpha$.

By lemma~\ref{graphegene}, there exist $m$ and $n$ such that the matrix $A^{\otimes m}\otimes B\otimes A^{\otimes n}$ has only finite entries and its critical graph has only one node.
By lemma~\ref{puissrg1}, there exists a neighborhood $V$ of $A^{\otimes m}\otimes B\otimes A^{\otimes n}$ and an integer $N$ such that every matrix in $V^{\otimes N}$ has rank $1$. Let $V_1\times V_2$ be a neighborhood of $(A,B)$ such that $V_1^{\otimes m}\otimes V_2\otimes V_1^{\otimes n}\subset V$. The matrices in $\left(V_1^{\otimes m}\otimes V_2\otimes V_1^{\otimes n}\right)^{\otimes N}$ have rank 1.
Since $A$ and $B$ are in the support of $\mu$, we have:
$$\P\left(rk\left(A\left((n+m+1)N\right)\otimes\cdots\otimes A(1)\right)=1\right)\ge \P\left(\left(A\left((n+m+1)N\right),\cdots, A(1)\right)\in \left(V_1^{m}\times V_2\times V_1^{ n}\right)^{ N}\right)>0.$$
\end{proofof}

\section{Powers of one matrix}
\subsection{General theory}\label{deter}
In this section, we briefly review the spectral and asymptotic theory of matrices with entries in $\rma$. For a complete exposition in English, see Baccelli et al.~\cite{BCOQ}.

\begin{definition}
A circuit on a directed graph is a closed path on the graph.
Let $A$ be a square matrix of size $k$ with entries in $\rma$.
\begin{enumerate}[i)]
\item The graph  of $A$ is the directed weighted graph whose nodes are the elements of $[1,k]$ and whose arcs are the $(i,j)$ such that $A_{ij}>-\infty$. The weight on $(i,j)$ is $A_{ij}$. The graph will be denoted by $\mathcal{G}(A)$ and the set of its elementary circuits  by $\mathcal{C}(A)$.
\item The weight of path $pth=(i_1,\cdots,i_n,i_{n+1})$ is $w(A,pth):=\sum_{j=1}^n A_{i_ji_{j+1}}.$ Its length is $|pth|:=n.$
\item The average weight of a circuit $c$ is $aw(A,c):=\frac{w(A,c)}{|c|}.$
\item The \mp spectral radius of $A$ is $\rhom(A):=\max_{c\in\mathcal{C}(A)}aw(A,c)$.
\item The critical graph of $A$ is obtained from $\mathcal{G}(A)$ by keeping only nodes and arcs belonging to circuits with average weight $\rhom(A)$. It will be denoted by $\mathcal{G}^c(A)$.
\item The cyclicity of a graph is the greatest common divisor of the length of its circuits if it is strongly connected (that is if every node can be reached from every other). Otherwise it is the least common multiple  of the cyclicities of its strongly connected components. The cyclicity of $A$ is that of $\mathcal{G}^c(A)$ and is denoted by $c(A)$.
\item The type of $A$ is \textbf{scsN-cycC}, where \textbf{N} is the number of strongly connected components of $\mathcal{G}^c(A)$ and \textbf{C} the cyclicity of $A$.
\end{enumerate}
\end{definition}

\begin{remark}\label{powerint} Interpretation of powers with $\mathcal{G}(A)$.\\
If $(i_1,i_2,\cdots,i_{n+1})$ is a path on $\mathcal{G}(A)$, its weight is $\sum_{j\in [1,n]}A_{i_ji_{j+1}}$, therefore $\left(A^{\otimes n}\right)_{ij}$ is the maximum  of the weights of paths from~$i$ to~$j$ with length~$n$.

Since the average weight of a circuit is an affine combination of the average weights of its minimal sub-circuits, the \mp spectral radius is the maximum of the average weights of all circuits.
\end{remark}

We first remind the (max,+)-spectral theory: if $\lambda\in\rma$ and $V\in \rma^k\backslash \{(-\infty)^k\}$ satisfy the equation $A\otimes V=\lambda\otimes V$, we say that $\lambda$ is an eigenvalue of $A$ and $V$ an eigenvector.

For every $A\in\M_k$, $\tilde{A}$ defined by $\tilde{A}_{ij}=A_{ij}-\rhom(A)$ satisfies $\rhom\left(\tilde{A}\right)=0$ and for every vector~$V$, we have $A\otimes V=\rhom(A)\otimes \tilde{A}\otimes V$, so it is enough to deal with the case $\rhom(A)=0$.  
To state the \mp spectral theorem, we need the following.

For every $A\in\rma^{k\times k}$ with $\rhom(A)\le 0$, we set:
$$A^+:=\bigoplus_{n\in [1,k]} A^{\otimes n}.$$

\begin{remark} For every $(i,j)\in [1,k]^2$, $A^+_{ij}$ is the maximum of the weights of paths from $i$ to $j$. Indeed, since $\rhom(A)\le 0$, all circuits have non-positive weights and removing circuits from a path makes its weight greater, so $\bigoplus_{n\in [1,k]} A^{\otimes n}  =\bigoplus_{n\ge 1} A^{\otimes n}$ and the remark follows from remark~\ref{powerint} .
\end{remark}

\begin{proposition}\label{propspec}Eigenvectors, G.~Cohen et al.~\cite{cohen83,cohen85a}.
\begin{enumerate}[i)]
\item If $c$ is a circuit on $\mathcal{G}^c(A)$, then its average weight is $\rhom(A)$.
\item If $\mathcal{G}(A)$ is strongly connected, then $\rhom(A)$ is the only eigenvalue of $A$.
\item If $\rhom(A)=0$, then for every $i\in\mathcal{G}^c(A)$, $A^+_{.i}$ is an eigenvector of $A$ with \mbox{eigenvalue $0$.}
\item If $\rhom(A)=0$, then for every eigenvector $y$ of $A$ with eigenvalue
$0$, we have $y=\bigoplus_{i\in\mathcal{G}^c(A)}y_i\otimes A^+_{.i}.$
\item If $\rhom(A)=0$, and if $i$ and $j$ are in the same strongly
connected component of $\mathcal{G}^c(A)$, then the column vectors $A^+_{.i}$ and $A^+_{.j}$
are proportional in the \mp sense.
\item If $\rhom(A)=0$, then no column vector $A^+_{.i}$ with
$i\in\mathcal{G}^c(A)$ is a max-plus linear combination of the
$A^+_{.j}$ with $j$ in other strongly connected components.
\end{enumerate}
\end{proposition}

\begin{proposition}\label{proppuiss}Powers, G.~Cohen et al.~\cite{cohen83,cohen85a}\\
Assume $\mathcal{G}(A)$ is strongly connected, $\rhom(A)=0$ and $c(A)=1$. Then, there is an $N\in\NN$ such that, for every $n\ge N$, we have $A^{\otimes n}=Q$,
where $Q$ is defined by
$\forall (i,j)\in
[1,k]^2,Q_{ij}:=\bigoplus_{l\in\mathcal{G}^c(A)}A^+_{il}\otimes
A^+_{lj}.$
\end{proposition}

\begin{remark}
For every $(i,j)\in [1,k]^2$, $Q_{ij}$ is the maximum of the weight of
paths from $i$ to $j$ that crosses~$\mathcal{G}^c(A)$.
\end{remark}

Lemma~\ref{scs1-cyc1} easily follows from propositions~\ref{propspec} and~\ref{proppuiss}:
\begin{proofof}{lemma~\ref{scs1-cyc1}} Let $\tilde{A}$ be the matrix defined by $\tilde{A}_{ij}=A_{ij}-\rhom(A)$. Then, $\rhom(\tilde{A})=0$ and $c(\tilde{A})=c(A)=1$. By proposition~\ref{proppuiss}, when $n$ is greater than some $N\in\NN$, the column vectors of $\tilde{A}^{\otimes n}$ are \mp eigenvectors of $\tilde{A}$. Therefore, by proposition~\ref{propspec}, all these vectors are proportional in the \mp sense, so $\tilde{A}^{\otimes n}$ has rank~1. But for every $i,j\in [1,k]$, $\left(A^{\otimes n}\right)_{ij}=\left(\tilde{A}^{\otimes n}\right)_{ij}+n\rhom(A)$, therefore $A^{\otimes n}$ also has rank~1.
\end{proofof}

We end this section with a simple but crucial lemma:
\begin{lemma}\label{lemGcA}
For every $A\in \mP_k$, there is an $N\in\NN$ such that every path $pth$ from $i$ to $j$ on $\G(A)$ with length $n\ge N$, and weight $w(A,pth)=A^{\otimes n}_{ij}$ crosses $\Gc(A)$.
\end{lemma}
Since this lemma is implicit in the published proofs of proposition~\ref{proppuiss} (G.~Cohen et al.~\cite{cohen85a}, F.~Baccelli et al.~\cite{BCOQ}), we prove it for sake of completeness:

\begin{proof}
Without loss of generality, we assume $\rhom(A)=0$.
By the definition of $\mP_k$, there is $N_1\in\NN$ such that $A^{\otimes N_1}$ has only finite entries. Let $M_1$ be the minimum of those entries.

Let $c$ be a circuit on $\Gc(A)$, with length $|c|$. According to proposition~\ref{propspec}, there exists an index $l\in [1,k]$ such that $A^{\otimes |c|}_{ll}=|c|\rhom(A)=0$. For every $n\in\NN$, and every $i,j\in [1,k]$, we have:
$$A^{\otimes (2N_1+n|c|)}_{i,j}\ge A^{\otimes N_1}_{il}+nA^{\otimes |c|}_{ll}+A^{\otimes N_1}_{lj}\ge 2 M_1.$$

Let $-\epsilon$ be the maximal average weight of circuits on $\G(A)$ not crossing $\Gc(A)$.
Let $pth=(i_1,\cdots,i_{n})$  be a path on $\G(A)$ not crossing $\Gc(A)$. It splits into a path of length at most $k$ and elementary circuits, with average weight at most $-\epsilon$. Denoting by $M_2$ the greatest entry of $A$, we have:
$$w(A, pth)\le k|M_2| -\epsilon (n-k).$$
Every $N$ large enough to satisfy $k(|M_2|+\epsilon)-N\epsilon<2 M_1$ thus satisfies the conclusion of the lemma.
\end{proof}

\subsection{Almost powers}
This section will be devoted to the proof of lemma~\ref{puissrg1}. This proof is based the ideas from the proofs of propositions~\ref{propspec}~and~\ref{proppuiss}. 

To understand the powers of $A$, we considered their entries as the weights of paths on $\mathcal{G}(A)$, as explained in remark~\ref{powerint}. We want to do the same for products of several matrices, which means the arcs weights can be different at each step.

From now on, $\mathcal{G}$ will be the complete directed graph with the elements of $[1,k]$ as nodes. For every finite sequence of matrices $(^i\hspace{-4pt}A)_{i\in [1,n]}$ we set the following notations:
\begin{enumerate}[-]
\item The weight of a path $pth=(i_j)_{j\in [1,n+1]}$ on $\mathcal{G}$ (with respect to $(^i\hspace{-4pt}A)_{i\in [1,n]}$) is 
$w (pth):=\sum_{j\in [1,n]} {}^j\hspace{-4pt}A_{i_ji_{j+1}}$
\item A path is maximizing if its weight is maximal among the weights of paths with the same origin, the same end, and the same length.
\end{enumerate} 
With these definitions $(i_j)_{j\in [1,n+1]}$ is maximizing if and only if its weight is $\left({^1\hspace{-4pt}A}\otimes\cdots\otimes{^n\hspace{-4pt}A}\right)_{i_1i_{n+1}}$.\\

\begin{proofof}{lemma~\ref{puissrg1}}
Since $\mathcal{G}^c(A)$ has only one node, there exists an $l\in [1,k]$
such that: $\forall c\in \mathcal{C}(A)\backslash\{(l,l)\},A_{ll}>aw(A,c).$
Thus, there exists an $\epsilon>0$ such that  
$$\forall c\in \mathcal{C}(A)\backslash\{(l,l)\},A_{ll}-aw(A,c)>3\epsilon.$$
Let $V$ be the open ball with center $A$ and radius $\epsilon$ for infinity norm and $M$ be the maximum of the infinity norm on $V$.

Let us notice that every matrix $B\in V$ has the same critical graph as $A$. Let  $\tilde{B}$ be the matrix with \mp spectral radius $0$ defined by $\tilde{B}_{ij}=B_{ij}-B_{ll}$. Then, $\|\tilde{B}-\tilde{A}\|_\infty<2\epsilon$ and for any elementary circuit $c\neq(l,l)$, the average weight of $c$ satisfy $aw(\tilde{A},c)<-3\epsilon$

From now on, $(^i\hspace{-4pt}A)_{i\in [1,n]}$ will be in $V^n$ and the weights of paths will always be with respect to $(\tilde{^i\hspace{-4pt}A})_{i\in [1,n]}$.

Let $pth=(i_j)_{j\in [1,n+1]}$ be a path of length $n$ that does not cross $l$. It can  split into a path of length less than $k$ and  elementary circuits.
Since an elementary circuit $c$ that is not $(l,l)$ and has weight $ w (c)\le  aw(\tilde{A},c)|c|+2|c|\epsilon<-|c|\epsilon$, we have:
$$ w (pth)<-(n-k)\epsilon +2kM.$$
But for every $i,j\in [1,k]$,
$$\left(\tilde{^1\hspace{-4pt}A}\otimes\cdots\otimes\tilde{^n\hspace{-4pt}A}\right)_{ij}\ge  w ((i,l,\cdots,l,j))>-2M,$$
so there exists an $N$ such that every maximizing path of length $n\ge N$ crosses~$l$.

Let $pth=(i_j)_{j\in [1,n+1]}$ be a maximizing path  of length $n\ge 2N+1$. Since $(i_j)_{j\in [1,N+1]}$ is also maximizing, there is a $j_0\le N$ such that $i_{j_0}=l$. Since $
(i_j)_{n-N\le j\le n+1}$ is maximizing for $(^j\hspace{-4pt}A)_{n-N\le j\le n}$, there exists $n-N\le j_1\le n+1$ such that $i_{j_1}=l$. The path $
(i_j)_{j_0\le j\le j_1}$ is a circuit, so it can  split into elementary circuits. Since elementary circuits have a negative weight, except for $(l,l)$, the only sub-circuit of $(i_j)_{j_0\le j\le j_1} $ is $(l,l)$. Consequently for every $j$ between $j_0$ and $j_1$,  and particular for $N+1\le j\le n-N$, $i_j=l$, therefore
$$ w (pth)= w \left((i_j)_{1\le j\le N+1}\right)+ w \left((i_j)_{n-N\le j\le n+1}\right).$$
This means that:
$$\forall n\ge 2N+1,\forall i,j\in [1,k],  \left(\tilde{^1\hspace{-4pt}A}\otimes\cdots\otimes\tilde{^n\hspace{-4pt}A}\right)_{i,j}= \left(\tilde{^1\hspace{-4pt}A}\otimes\cdots\otimes\tilde{^N\hspace{-4pt}A}\right)_{il}+\left(^{n-N}\tilde{\hspace{-4pt}A}\otimes\cdots\otimes\tilde{^n\hspace{-4pt}A}\right)_{lj},$$
therefore  $rk\left(\tilde{^1\hspace{-4pt}A}\otimes\cdots\otimes\tilde{^n\hspace{-4pt}A}\right)=1$ and also
$rk\left(^1\hspace{-4pt}A\otimes\cdots\otimes ^n\hspace{-6pt}A\right)=1$.
\end{proofof}

\section{Proof of the main lemma}\label{seclem}
\subsection{Reduced matrices}
To prove lemma~\ref{graphegene}, we set the following notions of reduced matrices.
\begin{definition}
A matrix is called reduced if it has only non-positive entries and at least one zero on each row.
It is called strictly reduced if it has only non-positive entries and exactly one zero on each row.
\end{definition}

\begin{lemma}\label{matred}\ 
\begin{enumerate}[i)]
\item The set of reduced matrices is a semigroup. So is the set of strictly reduced matrices.
\item Every reduced matrix $A$ has \mp spectral radius $0$ and $\Gc(A)$ is made of the circuits of $\G(A)$ whose arcs have weight~$0$.
\end{enumerate}
\end{lemma}

\begin{proof}
$i)$ is obvious. If $A$ is reduced, its coefficients are non-positive, and so is $\rhom(A)$. 
It is possible to build by induction a sequence $(i_j)_{j\in\NN}\in [1,k]^\NN$ such that for every $j\in\NN$, $A_{i_ji_{j+1}}=0$. This sequence takes twice the same value,  let us say in $j_1$ and $j_2$. Therefore, $c=(i_j)_{j_1\le j\le j_2}$ is a circuit of $\G(A)$ with arcs of weight~$0$. In particular, $w(A,c)=0$, so $\rhom(A)\ge 0$.

This shows that $\rhom(A)=0$, and, since the entries of $A$ are non-positive, the last statement is obvious.
\end{proof}

\begin{definition}
For every $A\in\M_k$ such that $\G(A)$ is connected, we define:
\begin{enumerate}
\item A matrix $\tilde{A}$ defined by $\forall i,j\in [1,k], \tilde{A}_{ij}=A_{ij}-\rhom(A)$.
\item A circuit $c_A$ on $\G(A)$: the smallest elementary circuit of $\Gc(A)$ for the lexicographical order.
\item An integer $\kappa(A)\in [1,k]$: the smallest node in $c_A$.
\item A matrix $\bar{A}$ defined by:
\begin{equation}\label{eqdefAbar}
\forall i,j\in [1,k], \bar{A}_{ij}=\tilde{A}_{ij}-\tilde{A}^+_{i\kappa(A)}+\tilde{A}^+_{j\kappa(A)}.
\end{equation}
\end{enumerate}
\end{definition}

We define the hyperplanes of lemma~\ref{graphegene} by linear forms:
\begin{definition}
Let $\G$ be the  complete directed graph whose nodes are the elements of $[1,k]$.
\begin{enumerate}
\item $\mE_1$ is the set of linear forms $aw(.,c_1)-aw(.,c_2)$, where $c_1$ and $c_2$ are two elementary circuits of~$\G$ that are different when seen as directed graphs.
\item $\mE_2$ is the set of linear forms $$w(.,pth_1)-|pth_1|aw(.,c)-w(.,pth_2)+|pth_2|aw(.,c),$$
where $pth_1$ and $pth_2$ are two elementary circuits of $\G$ with the same initial node~$i$ and same final node~$\kappa\neq i$, and where $c$ is an elementary circuit of $\G$ that goes through $\kappa$.
\end{enumerate}
\end{definition}

\begin{lemma}\label{A->Abar}\ 
\begin{enumerate}
\item For every matrix $A\in\M_k$ with strongly connected $\G(A)$, $\bar{A}$ is reduced. In particular, the weights of the arcs of $\Gc(\bar{A})$ are~$0$.
\item If no linear form of $\mE_1\cup\mE_2$ vanishes at $A$, then $\bar{A}$ is strictly  reduced.
\item The $0$-form is not in $\mE_1\cup\mE_2$.
\end{enumerate}
\end{lemma}

\begin{proofof}{lemma~\ref{A->Abar}}\ 
\begin{enumerate}
\item $\G(\tilde{A})$ is strictly connected because so is $\G(A)$. Moreover $\Gc(\tilde{A})$ is the same non-weighted graph as $\Gc(A)$, therefore $\kappa(A)\in\Gc(\tilde{A})$. Eventually, $\rhom(\tilde{A})=0$, and, by proposition~\ref{propspec}, the column vector $\tilde{A}^+_{.\kappa(A)}$ is an eigenvector of $\tilde{A}$ with eigenvalue 0. This is equivalent to each of the following equations systems:\\
%{\setlength{\tablecolsep}{0pt}
\begin{tabular}{lllcl}
$\forall i\in [1,k],$ &&$\tilde{A}_{i\kappa(A)}$&=&$\max_j (\tilde{A}_{ij}+\tilde{A}^+_{j\kappa(A)}),$\\

\multirow{2}{2cm}{$\forall i\in [1,k],\Big\{$}& 
$\forall j\in [1,k],$& $\tilde{A}_{i\kappa(A)}$&$\ge$&$\tilde{A}_{ij}+\tilde{A}^+_{j\kappa(A)}$\\
&$\exists j\in [1,k],$& $\tilde{A}_{i\kappa(A)}$&=&$\tilde{A}_{ij}+\tilde{A}^+_{j\kappa(A)}$\\

\multirow{2}{2cm}{$\forall i\in [1,k],\Big\{$}& 
$\forall j\in [1,k],$& $\bar{A}_{ij}$&$\le$&$0$\\
&$\exists j\in [1,k],$& $\bar{A}_{ij}$&=&$0$.\\
\end{tabular}%}

The last system exactly means that $\bar{A}$ is reduced.

\item Let us assume that no form in $\mE_1\cup\mE_2$ vanishes at $A$.

For every $i\in [1,k]$, there exists a path $pth_1=(i_1,\cdots, i_{|pth_1|+1})$ on $\G(\tilde{A})$ from $i$ to $\kappa(A)$ with weight $\tilde{A}^+_{i\kappa(A)}$ and with a minimal length among all paths with those properties.

We will show that $j=i_2$ is the only solution of equation $\bar{A}_{ij}=0$, therefore $\bar{A}$ is strictly reduced. This equation is equivalent to:
\begin{equation}\label{eqdefj}
\tilde{A}^+_{i\kappa(A)}=\tilde{A}_{ij}+\tilde{A}^+_{j\kappa(A)}.
\end{equation}

Let $j\in [1,k]$ be a solution of this equation and let $pth=(j_1,\cdots, j_{|pth|+1})$ be a path from $j$ to $\kappa(A)$ with weight $\tilde{A}^+_{j\kappa(A)}$ and  with a minimal length among all paths with those properties.

Since $\rhom(\tilde{A})=0$, the circuits of $\G(\tilde{A})$ have a non-positive weight, therefore the minimality of the lengths implies that $pth_1$ and $pth$ are elementary.\\

If $i=\kappa(A)$, then $w(\tilde{A},pth_1)=w(\tilde{A},(i,pth))=\tilde{A}^+_{\kappa(A)\kappa(A)}=0$, therefore $(i,pth)$ and $pth_1$ are circuits on $\G^c(\tilde{A})$, and therefore also circuits on $\G^c(A)$. $(i,pth)$ can be split into elementary circuits of $\G^c(A)$. Let $pth_2$ be the first one. Then, $aw(A,pth_2)=aw(A,pth_1)=\rhom(A)$ and, since no linear form in $\mE_1$ vanishes at $A$, $pth_1=pth_2$, therefore $i_2=j$.\\

If $i\neq \kappa(A)$, then we set $pth_2=(i,pth)$, therefore $w(\tilde{A},pth_2)=\tilde{A}^+_{i\kappa(A)}=w(\tilde{A},pth_1)$ or, equivalently 
$$ w(A,pth_1)-|pth_1|\rhom(A)=w(A,pth_2)-|pth_2|\rhom(A),$$
or
$$w(A,pth_1)-|pth_1|aw(c_A)=w(A,pth_2)-|pth_2|aw(c_A).$$
Since no linear form of $\mE_2$ vanishes at $A$, $pth_1=pth_2$, and $i_2=j$, provided $pth_2$ is elementary.\\

Let us assume it is not. There exists $l\in [1,|pth|+1]$ such that $j_l=i$. Then, we have:
\begin{eqnarray*}
w(\tilde{A},(i,j_1,\cdots,j_l))&=&w(\tilde{A},(i,pth))-w(\tilde{A},(j_l,\cdots,j_{|pth|+1}))\\
&=&\tilde{A}^+_{i\kappa(A)}-w(\tilde{A},(i,j_{l+1},\cdots,j_{|pth|},\kappa(A)))\ge 0,
\end{eqnarray*}  
therefore $aw(\tilde{A},(i,j_1,\cdots,j_l))=\rhom(\tilde{A})=0$. Therefore, $(i,j_1,\cdots,j_l)$ is a circuit on $\Gc(A)$, so it can be split into elementary circuits on $\G^c(A)$. Let $c_1$ be one of these circuits. Since $aw(A,c_1)=\rhom(A)=aw(A,c_A)$ and  no linear form of $\mE_2$ vanishes at $A$, it proves $c_1=c_A$. Therefore, $\kappa(A)\in \{j_1,\cdots,j_l\}$. Since $j_l=i\neq \kappa(A)$, we have $j_{|pth|+1}=\kappa(A)\in \{j_1,\cdots j_{l-1}\}$, therefore $pth$ is not elementary. Since we already noticed that  $pth$ is elementary, the assumption that $pth_2$ is not elementary was false, and this concludes the proof of point $(ii)$.

\item Obviously, the zero-form is not an element of $\mE_1$. Let us prove it is not one of $\mE_2$ either.

Let $pth_1=(i_1,\cdots,i_{|pth_1|+1})$, $pth_2=(j_1,\cdots,j_{|pth_2|+1})$ and $c=(l_1,\cdots,l_{|c|},l_1)$ be elementary path such that $i_1=j_1=i$, $i_{|pth_1|+1}=j_{|pth_2|+1}=l_1$ and $i_1\neq l_1$.

Let us assume that $f=w(.,pth_1)-|pth_1|aw(.,c)-w(.,pth_2)+|pth_2|aw(.,c)=0$ and show that $pth_1=pth_2$.

Since $pth_1$ and $pth_2$ are elementary and not circuits, $l_1$ appears only once in each: as last node. So, there is no arc leaving $l_1$ on $pth_1$ or on $pth_2$ and the component in $A^\circ_{l_1,l_2}$ of $w(.,pth_1)$ and $w(.,pth_2)$ is zero. The component of $f$ is then $\frac{(|pth_2|-|pth_1|)}{|c|}A^\circ_{l_1l_2}$. This is zero, thus $|pth_2|=|pth_1|$ and $f=w(.,pth_1)-w(.,pth_2)$.

Since $pth_1$ is elementary, $(i_1,i_2)$ is the only arc of $pth_1$ leaving $i$. Therefore, the component of $w(.,pth_1)$ in $A^\circ_{ii_2}$ is $A^\circ_{ii_2}$ and for every $j\neq i_2$, the component of $w(.,pth_1)$ in $A^\circ_{ij}$ is zero.

For the same reason, the component  of $w(.,pth_2)$ in $A^\circ_{ij_2}$ is $A^\circ_{ij_2}$ and for every  $j\neq j_2$, the component of $w(.,pth_2)$ in $A^\circ_{ij}$ is zero.
 
Since $w(.,pth_2)=w(.,pth_1)$, it follows from the last statements that $i_2=j_2$ and $w(.,(i_2,\cdots,i_{|pth_1|+1}))=w(.,(j_2,\cdots,j_{|pth_1|+1}))$.

By a finite induction, $pth_1=pth_2$.
\end{enumerate}
\end{proofof}
\subsection{Matrix with dominating diagonal}
Let $A\in\mP_k$ and $B\in\M_k$ be two matrices. By lemma~\ref{A->Abar}, we associate to $A$ a reduced matrix $\bar{A}$ defined by equation~(\ref{eqdefAbar}).
We also define $\hat{B}$ by:
\begin{equation}\label{eqdefBbar}\index{$\hat{B}$}
\hat{B}_{ij}:=B_{ij} -\tilde{A}^+_{i\kappa(A)}+\tilde{A}^+_{j\kappa(A)}.
\end{equation}
Warning: $\hat{B}$ also depends on $A$.

We want to apply the next lemma with $\bar{A}$ and $\hat{B}$. We will show that if some linear forms do not vanish at $(A,B)$, then $\bar{A}$ and $\hat{B}$ satisfy the hypotheses of this lemma.
\begin{lemma}\label{lemGcM}
Let $A\in\mP_k$ be a reduced matrix such that $\Gc(A)$ is strongly connected. Let $N$ be the integer given by lemma~\ref{lemGcA}\,.

Let $B\in\M_k$ be a matrix such that
\begin{equation}\label{eqANB}
\forall l_1,l_2,m_1,m_2\in [1,k], \left(A^{\otimes N} \otimes B\right)_{l_1m_1}=\left(A^{\otimes N} \otimes B\right)_{l_2m_2}\Rightarrow m_1=m_2.
\end{equation}

Then, there is $s\in\NN$ such that for every $p\in s+c(A)\NN$, $\Gc\left(A^{\otimes N} \otimes B\otimes A^{\otimes p}\right)$ is a complete directed graph.
Especially, $A^{\otimes N} \otimes B\otimes A^{\otimes p}$ has type~\textbf{scs1cyc1}.

Moreover, if $A$ is strictly reduced, then $\Gc\left(A^{\otimes N} \otimes B\otimes A^{\otimes p}\right)$ has exactly one node.
\end{lemma}

\begin{proof}
We first study the maximal entries of $ A^{\otimes N} \otimes B\otimes   A^{\otimes p}$.
Such an entry $\left( A^{\otimes N} \otimes B\otimes   A^{\otimes p} \right)_{ij}$ is the weight of a path $(i_r)_{r\in [1,N+p+2]}$ from $i$ to $j$ (figure~\ref{cheminmax}). Let $l$ be $i_{N+1}$ and  $m$ be $i_{N+2}$.

\begin{figure}[htb]
\begin{center}
\caption{Maximal paths of ${A}^{\otimes N}\otimes{B}\otimes {A}^{\otimes p}$.}\label{cheminmax}
\begin{picture}(0,0)%
\includegraphics{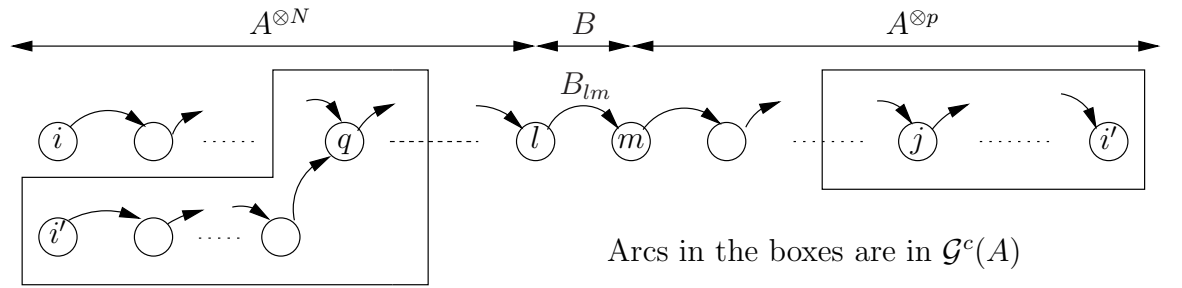}%
\end{picture}%
\setlength{\unitlength}{3947sp}%
\begingroup\makeatletter\ifx\SetFigFont\undefined%
\gdef\SetFigFont#1#2#3#4#5{%
  \reset@font\fontsize{#1}{#2pt}%
  \fontfamily{#3}\fontseries{#4}\fontshape{#5}%
  \selectfont}%
\fi\endgroup%
\begin{picture}(7149,1827)(289,-1873)
\put(7201,-1011){\makebox(0,0)[b]{\smash{\SetFigFont{12}{14.4}{\rmdefault}{\mddefault}{\updefault}{\color[rgb]{0,0,0}$i'$}%
}}}
\put(4201,-1001){\makebox(0,0)[b]{\smash{\SetFigFont{12}{14.4}{\rmdefault}{\mddefault}{\updefault}{\color[rgb]{0,0,0}$m$}%
}}}
\put(6001,-1001){\makebox(0,0)[b]{\smash{\SetFigFont{12}{14.4}{\rmdefault}{\mddefault}{\updefault}{\color[rgb]{0,0,0}$j$}%
}}}
\put(3601,-1001){\makebox(0,0)[b]{\smash{\SetFigFont{12}{14.4}{\rmdefault}{\mddefault}{\updefault}{\color[rgb]{0,0,0}$l$}%
}}}
\put(2401,-1001){\makebox(0,0)[b]{\smash{\SetFigFont{12}{14.4}{\rmdefault}{\mddefault}{\updefault}{\color[rgb]{0,0,0}$q$}%
}}}
\put(601,-1001){\makebox(0,0)[b]{\smash{\SetFigFont{12}{14.4}{\rmdefault}{\mddefault}{\updefault}{\color[rgb]{0,0,0}$i$}%
}}}
\put(601,-1601){\makebox(0,0)[b]{\smash{\SetFigFont{12}{14.4}{\rmdefault}{\mddefault}{\updefault}{\color[rgb]{0,0,0}$i'$}%
}}}
\put(4051,-1711){\makebox(0,0)[lb]{\smash{\SetFigFont{12}{14.4}{\rmdefault}{\mddefault}{\updefault}{\color[rgb]{0,0,0}Arcs in the boxes are in $\mathcal{G}^c(A)$}%
}}}
\put(1801,-286){\makebox(0,0)[lb]{\smash{\SetFigFont{12}{14.4}{\rmdefault}{\mddefault}{\updefault}{\color[rgb]{0,0,0}${A}^{\otimes N}$}%
}}}
\put(3826,-286){\makebox(0,0)[lb]{\smash{\SetFigFont{12}{14.4}{\rmdefault}{\mddefault}{\updefault}${B}$}}}
%\put(4951,-286){\makebox(0,0)[lb]{\smash{\SetFigFont{12}{14.4}{\rmdefault}{\mddefault}{\updefault}{\color[rgb]{0,0,0}${A}^{\otimes k}$}%
%}}}
\put(5801,-286){\makebox(0,0)[lb]{\smash{\SetFigFont{12}{14.4}{\rmdefault}{\mddefault}{\updefault}{\color[rgb]{0,0,0}${A}^{\otimes p}$}%
}}}
\put(3751,-661){\makebox(0,0)[lb]{\smash{\SetFigFont{12}{14.4}{\rmdefault}{\mddefault}{\updefault}${B}_{lm}$}}}
\end{picture}
\end{center}\end{figure}
Since $A$ is reduced, we define step by step a sequence of indices $\tilde{i}_r$ for $r$ between $N+2$ and $N+p+2$ such that $\tilde{i}_{N+2}=m$ and for every $r$, $A_{\tilde{i}_r\tilde{i}_{r+1}}=0$.
If we replace the $i_r$ with $r\ge N+3$ by $\tilde{i}_r$, we replace arcs with non positive weight by arcs with weight~0, so we get a path whose weight is greater or equal to $\left( A^{\otimes N} \otimes B\otimes   A^{\otimes p} \right)_{ij}$. Since this weight can not be strictly greater, and since $A$ is reduced, it means that $A_{i_ri_{r+1}}=0$. Eventually, the maximal entry $\left( A^{\otimes N} \otimes B\otimes   A^{\otimes p} \right)_{ij}$ is equal to $A^{\otimes N}_{il}+ B_{lm}$ and its value does not depend on $p$.

The choice of $N$ ensures that $(i_r)_{r\in [1,N+1]}$ crosses $\mathcal{G}^c(A)$. Let $q$ be the first node of the path on $\mathcal{G}^c(A)$. Starting from $q$ and following backward a circuit on $\Gc(A)$, we define a new path such that all nodes before $q$ are in $\mathcal{G}^c(A)$. Let~$i'$ be the starting point of this path (figure~\ref{cheminmax}). By construction, $i'\in \mathcal{G}^c(A)$ and $\left( A^{\otimes N}  \otimes B\right)_{i'm}$ is greater or equal to $\left( A^{\otimes N}  \otimes B   \right)_{im}$, therefore it is a maximal entry of~$A^{\otimes N}\otimes B$.\\

It follows from lemma~\ref{matred}~$ii)$ and the strong connectivity of $\mathcal{G}^c(A)$, that there exists $s_1\in N$ such that $A^{\otimes s_1}_{mi'}=0$. Since $c(A)$ is the cyclicity of the only strongly connected component of $\Gc(A)$, there exists $s_2\in\NN$, such that for every $p\in s_2+c(A)\NN$, we have $A^{\otimes p}_{i'i'}=0$. Therefore, setting $s=s_1+s_2$,  we have $A^{\otimes p}_{mi'}=0$ for every $p\in s+c(A)\NN$. From now on, we assume $p\in s+c(A)\NN$.

Matrix $M:= A^{\otimes N}  \otimes B\otimes   A^{\otimes p}$ has a maximal entry on its diagonal. This entry $M_{i'i'}$ is the  \mp spectral radius of $M$ and the weight of every arc of $\G^c(M)$. Moreover,  $M_{ij}=M_{i'i'}$ if and only if there  exists $m\in [1,k]$, such that $(A^{\otimes N}\otimes B)_{im}=M_{i'i'}$ and $A^{\otimes p}_{mj}=0$. Equation~(\ref{eqANB}) ensures that $m$ is the same for all~$i$.

If $i$ and $j$ are on $\Gc(M)$, then there exists an arc of $\Gc(M)$ from $i$ and another one from $j$, so that $(A^{\otimes N}\otimes B)_{im}=(A^{\otimes N}\otimes B)_{jm}=M_{i'i'}$. There are also arcs of $\Gc(M)$ to $i$ and to $j$, so that $A^{\otimes p}_{mi}=A^{\otimes p}_{mj}=0$. Therefore, $M_{ij}=M_{ji}=M_{i'i'}$ and $(i,j)$ and $(j,i)$ are arcs of~$\Gc(M)$. Eventually, $\Gc(M)$ is the complete directed graph whose nodes are the $i\in [1,k]$, such that $(A^{\otimes N}\otimes B)_{im}=M_{i'i'}$ and $A^{\otimes p}_{mi}=0$.

If $A$ is strictly reduced, so is $A^{\otimes p}$. Therefore, there is only one $i$ such that $A^{\otimes p}_{mi}=0$ and $\Gc(M)$ has only one node and one arc.
\end{proof}

To end the proof of lemma~\ref{graphegene}, we still have to show that if no linear form on $\mP_k\times\M_k$ in some finite set vanishes at $(A,B)$, then the matrices $\bar{A}$ and $\hat{B}$ defined by equations~(\ref{eqdefAbar}) and~(\ref{eqdefBbar}) satisfy the hypotheses of lemma~\ref{lemGcM}.
Actually, for every circuit $c$, we have~
$$aw\left(c,\bar{A}^{\otimes N} \otimes\hat{B}\otimes  \bar{A}^{\otimes p}\right)=aw\left(c,A^{\otimes N}  \otimes B\otimes   A^{\otimes p}\right)+(N+p)\rhom(A),$$
Therefore, $\bar{A}^{\otimes N} \otimes\hat{B}\otimes  \bar{A}^{\otimes p}$ and $A^{\otimes N}  \otimes B\otimes   A^{\otimes p}$ have the same critical graph, with different weights,  and for $p$ large enough, $A^{\otimes p}\in\R^{k\times k}$, therefore $A^{\otimes N}  \otimes B\otimes   A^{\otimes p}\in\R^{k\times k}$.

We set the following definition:
\begin{definition}
Let $\G$ be the complete directed graph with node in~$[1,k]$ and $A^\circ$ (resp. $B^\circ$) the function, that maps $(A,B)\in\mP_k\times\M_k$ to $A$ (resp. $B$).

We denote by $\mE_3$ the set of linear forms on $\mP_k\times\M_k$ of the following type:
\begin{equation}\label{eqflE3}
{\arraycolsep2pt\begin{array}{rl}
B^\circ_{j_1m_1}-B^\circ_{j_2m_2}
+w(A^\circ,pth_{i_1j_1})&-w(A^\circ,pth_{i_1\kappa})+w(A^\circ,pth_{m_1\kappa})\\
-w(A^\circ,pth_{i_2j_2})&+w(A^\circ,pth_{i_2\kappa})-w(A^\circ,pth_{m_2\kappa})\\
-(|pth_{i_1j_1}|-|pth_{i_1\kappa}|+|pth_{m_1\kappa}|-&|pth_{i_2j_2}|+|pth_{i_2\kappa}|-|pth_{m_2\kappa}|)aw(A^\circ,pth_{\kappa\kappa}),
\end{array}}
\end{equation}
where $i_1,i_2,j_1,j_2,m_1,m_2$ and $\kappa$ are nodes of $\G$ such that $m_1\neq m_2$ and for every $i\in\{i_1,i_2,m_1,m_2,\kappa\}$,
$pth_{i\kappa}$ is an elementary path on $\G$ from~$i$ to~$\kappa$ and for every $l\in [1,2]$, $pth_{i_lj_l}$ is a path on $\G$ from~$i_l$ to~$j_l$ with length at most~$|pth_{\kappa\kappa}|k$.
\end{definition}

Since $\mE_3$ obviously does not contain the zero-form, lemma~\ref{graphegene} follows from the next lemma:
\begin{lemma}\label{lemAbarBbar}\ 
\begin{enumerate}
\item If no form in $\mE_1$ vanishes at $A\in\mP_k$ , then $\Gc(\bar{A})$ is strongly connected.
\item If no form in $\mE_3$ vanishes at  $(A,B)\in\mP_k\times\M_k$, then $(\bar{A},\hat{B})$ satisfy relation~(\ref{eqANB}).
\end{enumerate}
\end{lemma}

\begin{proof}\ 
\begin{enumerate}
\item 
Since $\Gc(A)$ and $\Gc(\bar{A})$ are equal as non-weighted graphs, we only have to show that $\Gc(A)$ is strongly connected. But no form in  $\mE_1$ vanishes at $A$  means that every two elementary circuits of $\G(A)$ have distinct weights. Therefore, $\Gc(A)$ is an elementary circuit, so it is strongly connected.
\item Let $(A,B)\in\mP_k\times\M_k$ be such that $(\bar{A},\hat{B})$ does not satisfy relation~(\ref{eqANB}). Then, there exist $i_1,i_2,m_1,m_2\in [1,k]$ such that $m_1\neq m_2$ and:
$$\left(\bar{A}^{\otimes N} \otimes\hat{B}\right)_{i_1m_1}=\left(\bar{A}^{\otimes N} \otimes\hat{B}\right)_{i_2m_2}.$$

For each $l\in [1,2]$, we take $j_l\in [1,k]$ such that $\left(\bar{A}^{\otimes N} \otimes\hat{B}\right)_{i_lm_l}=\left(\bar{A}^{\otimes N}\right)_{i_lj_l}+\hat{B}_{j_lm_l}$ and we denote $\kappa(A)$ by $\kappa$.

For each $i\in\{i_1,i_2,m_1,m_2,\kappa\}$, we take  $pth_{i\kappa}$ a path on $\G(A)$ from~$i$ to~$\kappa$  such that $w(\tilde{A},pth_{i\kappa})=\tilde{A}^+_{i\kappa}$ with minimal length among such paths.

For each $l\in [1,2]$, we take  $pth_{i_lj_l}$ a path on  $\G(A)$ from~$i_l$ to~$j_l$  such that $w(\tilde{A},pth_{i_lj_l})=\tilde{A}^{\otimes N}_{i_lj_l}$ with minimal length among such paths.

Since $\rhom(\tilde{A})=0$, the circuits on $\G(A)$ have non-positive weight, so the minimality of length ensure that for every $i\in\{i_1,i_2,m_1,m_2,\kappa\}$, the path $pth_{i\kappa}$  is elementary.

Proposition~\ref{proppuiss} applied to $\tilde{A}^{\otimes c(A)}$ ensures that $|pth_{i_lj_l}|\le 2kc(A)$. And $c(A)$ divides $|pth_{\kappa\kappa}|$, so $|pth_{i_lj_l}|\le 2k|pth_{\kappa\kappa}|$.

With those notations, it follows from the definition of $\bar{A}$ and $\hat{B}$ that the linear form $f\in\mE_3$ defined by formula~(\ref{eqflE3}) vanishes at $(A,B)$.
\end{enumerate}
\end{proof}

\begin{remark}
Putting all lemmas together, we get the following semi-explicit formulation of theorems~\ref{thgenedis} and~\ref{thgenesupp}:

Let $(A,B)\in \mP_k\times\M_k$ be such that no linear form in $\mE_1\cup\mE_2$ vanishes at $A$ and no linear form in $\mE_3$ vanishes at $(A,B)$.
Then, every sequence $\left(A(n)\right)_{n\in\NN}$ of i.i.d. random matrices taking values $A$ and $B$ with positive probability has the MLP.

Let $\mu$ be a probability measure on $\M_k$ with support $S_\mu$. If there is $A\in \mP_k\cap S_\mu$ at which no linear form in $\mE_1\cup\mE_2$ vanishes, and if $S_\mu$ is not included in a union of finitely many sets of type $\{A\in\M_k|A_{ij}-A_{\kappa l}=a\}$, with $i,j,\kappa,l\in [1,k]$, $j\neq l$ and $a\in\R$, then every sequence $\left(A(n)\right)_{n\in\NN}$ of i.i.d. random matrices with law $\mu$ has the MLP.
\end{remark}

\bibliographystyle{amsplain}
\bibliography{max+}

\def\cprime{$'$} \def\cprime{$'$}
\providecommand{\bysame}{\leavevmode\hbox to3em{\hrulefill}\thinspace}
\providecommand{\MR}{\relax\ifhmode\unskip\space\fi MR }
\providecommand{\MRhref}[2]{%
  \href{http://www.ams.org/mathscinet-getitem?mr=#1}{#2}
}
\providecommand{\href}[2]{#2}
\begin{thebibliography}{10}

\bibitem{Baccelli}
F.~Baccelli, \emph{Ergodic theory of stochastic {P}etri networks}, Ann. Probab.
  \textbf{20} (1992), no.~1, 375--396. \MR{MR1143426 (93a:68110)}

\bibitem{BCOQ}
F.~Baccelli, G.~Cohen, G.J. Olsder, and J.P. Quadrat, \emph{Synchronisation and
  linearity}, John Wiley and Sons, 1992.

\bibitem{braker}
H.~Braker, \emph{Algorithms and applications in timed discrete event systems},
  Ph.D. thesis, Delft University of Technology, Dec 1993.

\bibitem{cohen83}
G.~Cohen, D.~Dubois, J.P. Quadrat, and M.~Viot, \emph{Analyse du comportement
  p\'{e}riodique des syst\`{e}mes de production par la th\'{e}orie des dio\"\i
  des}, Rapport de recherche 191, INRIA, Le Chesnay, France, 1983.

\bibitem{cohen85a}
\bysame, \emph{A linear system theoretic view of discrete event processes and
  its use for performance evaluation in manufacturing}, IEEE Trans. on
  Automatic Control \textbf{AC--30} (1985), 210--220.

\bibitem{Cohen}
J.~E. Cohen, \emph{Subadditivity, generalized products of random matrices and
  operations research}, SIAM Rev. \textbf{30} (1988), no.~1, 69--86.
  \MR{MR931278 (89g:60300)}

\bibitem{GaubertHong}
S.~Gaubert and D.~Hong, \emph{Series expansions of lyapunov exponents and
  forgetful monoids}, Tech. report, INRIA, 2000.

\bibitem{GaubertMairesseIEEE}
S.~Gaubert and J.~Mairesse, \emph{Modeling and analysis of timed {P}etri nets
  using heaps of pieces}, IEEE Trans. Automat. Control \textbf{44} (1999),
  no.~4, 683--697. \MR{MR1684424 (99m:68139)}

\bibitem{gaumair95}
St{\'e}phane Gaubert and Jean Mairesse, \emph{Task resource models and
  {$(\max,+)$} automata}, Idempotency (Bristol, 1994), Publ. Newton Inst.,
  vol.~11, Cambridge Univ. Press, Cambridge, 1998, pp.~133--144. \MR{MR1608394
  (98j:68031)}

\bibitem{CaractMpQueuNet}
B.~Heidergott, \emph{A characterisation of {$(\max,+)$}-linear queueing
  systems}, Queueing Systems Theory Appl. \textbf{35} (2000), no.~1-4,
  237--262. \MR{MR1782609 (2001h:90014)}

\bibitem{HeidergottDeVriesPubTransNet}
B.~Heidergott and R.~de~Vries, \emph{Towards a ({M}ax,+) control theory for
  public transportation networks}, Discrete Event Dyn. Syst. \textbf{11}
  (2001), no.~4, 371--398. \MR{MR1852744 (2002f:93074)}

\bibitem{MpAtWork}
B.~Heidergott, G.~J. Oldser, and J.~van~der Woude, \emph{Max plus at work},
  Princeton Series in Applied Mathematics, Princeton University Press,
  Princeton, NJ, 2006, Modeling and analysis of synchronized systems: a course
  on max-plus algebra and its applications. \MR{MR2188299 (2006g:93079)}

\bibitem{Mairesse}
J.~Mairesse, \emph{Products of irreducible random matrices in the {$(\max,+)$}
  algebra}, Adv. in Appl. Probab. \textbf{29} (1997), no.~2, 444--477.
  \MR{MR1450939 (98k:60165)}

\bibitem{limtop}
G.~Merlet, \emph{Limit theorems for iterated random topical operators}, Tech.
  report, IRMAR, http://hal.ccsd.cnrs.fr/ccsd-00004594, 2005, Submitted.

\bibitem{TclGM}
\bysame, \emph{A central limit theorem for stochastic recursive sequences of
  topical operators}, Tech. report, Keio University,
  http://hal.ccsd.cnrs.fr/ccsd-00082113, 2006.

\bibitem{vincent}
J.-M. Vincent, \emph{Some ergodic results on stochastic iterative discrete
  events systems}, Discrete Event Dynamic Systems \textbf{7} (1997), no.~2,
  209--232.

\end{thebibliography}
\end{document}